\RequirePackage{fix-cm}
\documentclass[smallextended]{svjour3}       
\smartqed  
\usepackage{graphicx}
\usepackage{amsmath,amssymb,bbm,url}
\usepackage{pgf, tikz}
\usetikzlibrary{arrows, automata}

\spnewtheorem*{conj}{Conjecture}{\bf}{\it}

\newcommand{\reals}{\mathbb{R}}

\DeclareMathOperator{\Prox}{Prox}

\DeclareMathOperator*{\dom}{dom}
\DeclareMathOperator*{\dist}{dist}
\DeclareMathOperator*{\ran}{range}
\DeclareMathOperator*{\rec}{rec}

\DeclareMathOperator*{\ri}{ri}
\DeclareMathOperator*{\argmin}{arg\,min}
\newcommand{\vz}{0}

\begin{document}

\title{Douglas--Rachford Splitting and ADMM for Pathological Convex Optimization
}

\titlerunning{DRS and ADMM for Pathological Convex Optimization}        

\author{Ernest K. Ryu         \and
        Yanli Liu \and Wotao Yin 
}


\institute{Ernest K. Ryu \at
              UCLA\\
              \email{eryu@math.ucla.edu}           
           \and
Yanli Liu \at
              UCLA\\
              \email{yanli@math.ucla.edu}           
           \and
Wotao Yin \at
              UCLA\\
              \email{wotaoyin@math.ucla.edu}           
}

\date{Received: date / Accepted: date}

\makeatletter

\renewcommand{\paragraph}{%
  \@startsection{paragraph}{4} {\z@}
  {1ex \@plus 1ex \@minus .2ex}
  {-1ex}%
  {\bfseries}%
}                                  
 
\makeatother

\maketitle

\begin{abstract}
Despite the vast literature on DRS and ADMM, there has been very little work analyzing their behavior under pathologies. Most analyses assume a primal solution exists, a dual solution exists, and strong duality holds. When these assumptions are not met, i.e., under pathologies, the theory often breaks down and the empirical performance may degrade significantly. In this paper, we establish that DRS only requires strong duality to work, in the sense that asymptotically iterates are approximately feasible and approximately optimal.


\end{abstract}

\section{Introduction}
\label{s:intro}
Douglas--Rachford splitting (DRS)
and alternating directions method of multipliers (ADMM)
are classical methods originally presented in 
\cite{peaceman1955,douglas1956,lions1979,kellogg1969} 
and \cite{gabay1976,glowinski1975}, respectively.
DRS and ADMM are closely related.
Over the last decade, these methods
have enjoyed a resurgence of popularity,
as the demand to solve ever larger problems grew.

DRS and ADMM have strong theoretical guarantees and empirical performance,
but such results are often limited to non-pathological problems; 
in particular, most analyses assume
a primal solution exists, a dual solution exists, and strong duality holds.
When these assumptions are not met, i.e., under pathologies,
the theory often breaks down and the empirical performance may degrade significantly.
Surprisingly, there had been very little work 
analyzing DRS and ADMM under pathologies,
despite the vast literature on these methods.
There has been some recent exciting progress in this area, which we review in Section~\ref{ss:prior}.

In this paper, we analyze the asymptotic behavior of DRS and ADMM under pathologies.
While it is well known that the iterates ``diverge'' in such cases,
the precise manner in which they do so was not known.
We establish that
when strong duality holds, i.e., when
$p^\star=d^\star\in[-\infty,\infty]$,
DRS works,
in the sense that asymptotically the divergent iterates 
are approximately feasible and approximately optimal.
The assumption that 
primal and dual solutions exist is not necessary.
We then translate the pathological analyses for DRS
to pathological analyses for ADMM.

Furthermore, we conjecture that DRS necessarily fails when strong duality fails,
and we present empirical evidence that supports (but does not prove) this conjecture.
In other words, we believe
strong duality is the necessary and sufficient condition for DRS to work.

%
%
%
%

\subsection{Summary of results, contribution, and organization}
Sections~\ref{s:app-drs} and \ref{s:app-admm} present what we consider the fruits of this work, 
the convergence analyses of DRS and ADMM under various pathologies.
In fact, we suggest readers read Sections~\ref{s:app-drs} and \ref{s:app-admm}
before reading the theory of Section~\ref{s:theory}, 
as doing so will give a sense of direction.

We quickly illustrate, through examples, the kinds of results we show.
Precise definitions and statements are presented later.
We want DRS and ADMM to find a point that is \emph{approximately feasible} and, when applicable, \emph{approximately optimal}.
For example, if the primal problem is weakly infeasible,  we want the DRS iterates to satisfy
\[
x^{k+1}-x^{k+1/2}\rightarrow \vz
\]
and we show this as Theorem~\ref{thm:drs-f}.
As another example, if the primal problem is feasible but has no solution and $d^\star=p^\star>-\infty$,  we want the DRS iterates to satisfy
\begin{gather*}
x^{k+1}-x^{k+1/2}\rightarrow \vz,\qquad
f(x^{k+1/2})+g(x^{k+1})\rightarrow p^\star
\end{gather*}
and we show related results
as Theorems~\ref{thm:drs-c} and \ref{thm:drs-d}.
We can say something for all the pathological cases, so long as $d^\star=p^\star$.

Section~\ref{s:theory} presents the main theoretical contribution of this work.
To show that DRS and ADMM successfully achieve the 2 goals of
approximate feasibility and approximate optimality,
we need 2 separate major theoretical components.

Section~\ref{s:infim} presents the first component, which
analyzes the ``fixed-point iteration'' without a fixed point
with tools from operator theory.
With this machinery, we show results like
$x^{k+1}-x^k\rightarrow \vz$
or
$x^{k+1}-x^k\rightarrow v$, where $v$ is a certificate of (primal or dual) infeasibility.
Our contribution is defining the notion of improving directions via recession functions
and fully characterizing the infimal displacement vector with this notion.

Section~\ref{ss:f-val-anal} presents the second component,
the function-value analysis,
which is based on ideas from convex optimization and subgradient inequalities.
With these techniques, we show results like
$f(x^{k+1/2})+g(x^{k+1})\rightarrow p^\star$.
This part requires the $d^\star=p^\star$ assumption.
Our function-value analysis uses, but does not immediately follow from,
the results of Section~\ref{s:infim}.
To the best of our knowledge,
analyzing the convergence of objective values for DRS or ADMM applied to pathological problems
has not been done before.

Section~\ref{ss:shadow} presents a third, relatively minor theoretical component,
which we use later in Section~\ref{s:app-admm} to
translate analyses for DRS to analyses for ADMM.

As the goal of this work is to prove several theorems, one each for the many pathological cases,
we  build up our theory in a series of lemmas and corollaries.
Some of these lemmas are rather simple extensions of known results while some are novel.
All results of Section~\ref{s:theory} are eventually 
used in proving the
5 theorems of Section~\ref{s:app-drs}
and the 
3 theorems of Section~\ref{s:app-admm}.

The paper is organized as follows.
Section~\ref{s:prelim} reviews standard notions of convex analysis,
states several known results, and sets up the notation.
Section~\ref{s:theory} presents the main theoretical contribution of this paper.
Section~\ref{s:app-drs}
analyzes DRS under pathologies 
with the theory of Section~\ref{s:theory}.
Section~\ref{s:app-admm}
analyzes ADMM under pathologies 
with the theory of Sections~\ref{s:app-drs} and \ref{s:theory}.
Section~\ref{ss:sd-fail} presents counterexamples to make 
additional observations.
Section~\ref{s:conclusion} concludes the paper.

\subsection{Prior work}
\label{ss:prior}


As pathological convex optimization problems do arise in practice  \cite{lofberg2009,loera2012,drusvyatskiy2017,Waki2012,Waki20122},
there is practical value 
in studying how well-behaved and robust an algorithm is in such setups
However, the 
there had been surprisingly little work investigating the behavior of the popular methods DRS and ADMM under pathologies.
The understanding is still incomplete, but there has been some recent progress:
\cite{BauschkeCombettesLuke2004_finding,bauschke2016douglas,BauschkeMoursi2017_douglasrachford,liu2017new}
analyze DRS under specific pathological setups,
\cite{bauschke2014generalized,bauschke2016,BauschkeMoursi2017_douglasrachford}
analyze DRS under general setups,
and 
\cite{raghunathan2014,osqp,osqp-infeasibility} analyze ADMM
under specific pathological setups for conic programs.
These studies, however, are limited to more specific setups
and pathologies where an improving direction exists or the primal problem is strongly feasible.

The convex feasibility problem of finding an $x\in A\cap B$, where $A$ and $B$ are nonempty closed convex sets,
is a subclass of problems with practical importance.
While it is possible to recast convex feasibility problems into equivalent optimization problems and apply the results of this work, prior work on the specific setup has stronger results \cite{BauschkeCombettesLuke2004_finding,bauschke2016douglas,BauschkeMoursi2017_douglasrachford}.
We discuss further comparisons in Section~\ref{ss:feas}.

DRS has strong primal-dual symmetry,
in the sense of Fenchel duality
for convex optimization  \cite{fenchel1953,rockafellar1970} 
 and, more generally, 
 Attouch-Th\'era duality 
for monotone operators
\cite[p.\ 40]{mercier1980} and  \cite{attouch1996}.
See \cite[Lemma~3.6 p.\ 133]{eckstein1989} or   \cite{bauschke2012,Bauschke2017,BauschkeMoursi2017_douglasrachford}
for in-depth studies on this subject.
Naturally, our results also exhibit a degree of primal-dual symmetry,
although we do not explicitly address it in the interest of space.
Rather, we take the viewpoint that the primal problem is the problem of interest
and the dual problem is an auxiliary conceptual and computational tool.

In operator theory, 
and especially in infinite dimensional problems arising from physics and PDEs,
the sum of two maximal monotone operators may not be maximal,
and one can consider this a pathology.
One remedy to such pathology
is to generalize the notion of the sum of two operators
by regularizing the operators and then considering the limit as the regularization is reduced to zero
\cite{attouch1994,REVALSKI1999979,REVALSKI2002505}.
This notion of pathology and the remedy is quite different from what we consider.
For some of the pathologies we consider, $\partial f+\partial g$ is a perfectly well-defined
maximal monotone operator. Moreover, we do not remedy the pathology but rather simply analyze how DRS and ADMM
behave under the pathology.
We work in finite dimensions and thereby avoid the notion of weak and strong convergence.

When a problem is known to be pathological a priori,
one can first modify or regularize the problem and then solve the non-pathological problem.
One such approach is facial reduction, a pre-processing step that
rids a pathological conic program of difficult pathologies
\cite{borwein1981regularizing,borwein1981facial,borwein1982,ramana1997strong,pataki2000simple,cheung2013,waki2010,waki2013facial,permenter2014,lourencco2015solving,permenter2017solving,Permenter2017,pataki2017sieve}.
In contrast, the goal of this work is to analyze
DRS and ADMM when directly applied to pathological convex programs.
To put in differently, we do not assume
users of DRS or ADMM have a priori knowledge of whether the problem is pathological.


The standard analysis for DRS proves the iterates converge using ideas from operator theory and fixed point iterations
\cite{lions1979,eckstein1989,eckstein1992,combettes2004,combettes2011,eckstein2015,Combettes2018}.
The standard analyses of ADMM prove the iterates converge
by reducing ADMM to DRS \cite{gabay1983,eckstein1992,eckstein2015} or with a direct analysis via a Lyapunov function 
\cite{fortin1983,glowinski1984,bertsekaspar,Deng2016,Chen2017}.
These analyses rely on the existence of a primal-dual saddle point,
which only exists under the non-pathological case,
and therefore do not immediately generalize to pathological setups.

The first part of our analysis 
relies on a classical result by Pazy  \cite{Pazy1971_asymptotic} 
and Baillon et al.\ \cite{BaillonBruckReich1978_asymptotic} from the 1970s,
which characterize the asymptotic behavior of fixed-point iterations without fixed-points.
There has been some recent work that analyze algorithms that can be interpreted as fixed-point iterations without fixed-points
\cite{bauschke1997method,BauschkeCombettesLuke2004_finding,borwein2013cyclic,artacho2014douglas,bauschke2014generalized,bauschkepazy2016,bauschke2016douglas,artacho2016global,moursi2016forward,BauschkeMoursi2017_douglasrachford,ryu2017cosmic,liu2017new}. 
The analysis of Section~\ref{s:infim} was inspired by these works.

Another recent line of analysis for DRS and ADMM is function-value analysis,
which establishes the objective values, rather than the iterates, converge
\cite{davis2015,davis2016,davis2017}.
These analyses, however,  also rely on the existence of a primal-dual saddle point
and do not immediately generalize to the pathological setups.
The function-value analysis of Section~\ref{ss:f-val-anal} was inspired by these works.

\section{Preliminaries}
\label{s:prelim}
In this section, we review standard notions of convex analysis,
state several known results, and set up the notation.
For the sake of brevity,
we omit proofs or direct references of the standard results
and refer interested readers to standard references such as
\cite{rockafellar1970,rockafellar1974,BauschkeCombettes2017_convex}.

 Throughout this section, 
 $f:\reals^n\rightarrow \mathbb{R}\cup\{\infty\}$ is a function
 and $A,B,C\subseteq \reals^n$ are nonempty sets.
 
A set $C$ is convex, if
 $x,y\in C$ and $\theta\in[0,1]$ implies
$\theta x+(1-\theta)y\in C$.
Write $\overline{C}$ for the closure of $C$.
If $C$ is convex $\overline{C}$ is convex.
The Minkowski sum 
and differences of $A$ and $B$ are
\[
A+B=\{a+b\,|\,a\in A,\,b\in B\},
\qquad
A-B=\{a-b\,|\,a\in A,\,b\in B\},
\]
respectively. If $A$ and $B$ are convex, then
$A+B$ and $A-B$ are convex.
However, neither $A+B$ nor $A-B$ is guaranteed to be closed, even when $A$ and $B$ are nonempty closed convex sets.

For the distance between $x\in \reals^n$ and the set $A$,
write
\[
\dist(x,A)=\inf\{\|x-a\|\,|\,a\in A\}.
\]
For the distance between $A$ and $B$, write
\[
\dist(A,B)=\inf\{\|a-b\|\,|\,a\in A,\,b\in B\}.
\]
Note that $\dist(A,B)=0$ if and only if $\vz\in \overline{A-B}$.

A function $f$ is convex if
 $ f(\theta x+(1-\theta)y)\le \theta f(x)+(1-\theta)f(y)$
 for all $x,y\in \reals^n$ and $\theta\in[0,1]$.
A function $f$ is closed if its epigraph
$\big\{(x,\alpha)\in \reals^{n+1}\,\mid\,f(x)\le \alpha\big\}$
is a closed subset of $\mathbb{R}^{n+1}$.
We say $f:\reals^n\rightarrow \mathbb{R}\cup\{\infty\}$ is proper if $f(x)<\infty$ for some $x$.
In this paper, we focus our attention on closed, proper, and convex (CPC) functions.
 If $f$ and $g$ are CPC functions, then $f+g$ is CPC or $f+g=\infty$ everywhere.
If $\gamma>0$, then $\gamma f$ is CPC.

Define the (effective) domain of $f$ as
$ \dom f=\{x\in \mathbb{R}^n\,|\,f(x)<\infty\}$.
If $f$ is a convex function, then $ \dom f$ is a convex set.
However, $ \dom f$ may not be closed even when $f$ is CPC.
For any $\gamma>0$, we have $\dom \gamma f=\dom f$.

 Define the conjugate of $f$ as
 $ f^*(y)=\sup_{x\in \reals^n}\{ \langle y, x\rangle-f(x) \}$.
If $f$ is convex and proper,
  then $f^*:\reals^n\rightarrow \mathbb{R}\cup\{\infty\}$ is CPC.
If $f$ is CPC, then $(f^*)^*=f$.
For any $\gamma>0$, we have
$(\gamma f)^*(x)=\gamma f^*(x/\gamma)$
and $\dom(\gamma f)^*=\gamma \dom f^*$.
If $h(x)=g(-x)$, then $h^*(y)=g^*(-y)$.

Define the projection onto $C$ as 
$\Pi_C(x_0)=\argmin_{x\in C}\|x-x_0\|$.
When $C$ is closed and convex,  $\Pi_C:\reals^n\rightarrow\reals^n$ 
is well-defined, i.e., the minimizer uniquely exists.

Define the indicator function with respect to $C$ as
\[
\delta_C(x)=
\left\{
\begin{array}{ll}
0&\text{if } x\in C\\
\infty&\text{otherwise.}
\end{array}
\right.
\]
When $C$ is closed convex, 
$\delta_C:\mathbb{R}^n\rightarrow\reals\cup\{\infty\}$ is CPC.

Define the support function of $C$ as
\[
\sigma_C(y)=\sup_{x\in C}\{\langle x, y\rangle\}.
\]
$\sigma_C:\mathbb{R}^n\rightarrow\reals\cup\{\infty\}$ is CPC.
When $C$ is convex, we have 
$\sigma_C=\sigma_{\overline{C}}$.
If $A$ and $B$ are convex, 
then $\sigma_{A+B}=\sigma_A+\sigma_B$.
If $C$ is closed  and convex  then 
$(\sigma_C)^*=\delta_C$.

Define the recession function of $f$ as
 \begin{equation}
\rec f(d)=
\lim_{\alpha \rightarrow\infty} \frac{f(x+\alpha  d)-f(x)}{\alpha}
\label{eq:rec-def}
\end{equation}
for any $x\in \dom f$.
Loosely speaking, 
the recession function
characterizes the asymptotic change of $f$ as we go in direction $d$.
In fact,
\[
f(x+\alpha d)=\alpha \rec f(d)+o(\alpha)
\]
as $\alpha\rightarrow\infty$ for any $x\in \dom f$.
The recession function $\rec f: \reals^n\rightarrow \reals\cup\{\infty\}$ is a positively homogeneous CPC function.
If $h(x)=g(-x)$, then 
$\rec(h^*)(d)=\rec (g^*)(-d)$.
When $f$ and $g$ are CPC,
either $f(x)+g(x)=\infty$ for all $x\in \mathbb{R}^n$
or $\rec(f+g)=\rec f+\rec g$.
If $f$ is CPC,  then $\sigma_{\dom f^* }=\rec f$.
 Define the proximal operator 
 $\Prox_f:\reals^n\rightarrow\reals^n$ as
 \[
 \Prox_f(z)=\argmin_{x\in \reals^n}
 \left\{ f(x)+(1/2)\|x-z\|^2 \right\}.
 \]
When $f$ is CPC, the $\argmin$ uniquely exists, and therefore $\Prox_f$ is well-defined.
When $C$ is closed and convex, 
$ \Prox_{\delta_C}=\Pi_C$.
When $f$ is CPC, 
$\Prox_{f}+\Prox_{f^*}=I$, where $I:\reals^n\rightarrow\reals^n$ is the identity operator.

A mapping $T:\reals^n\rightarrow\reals^n$
is nonexpansive if
$\|T(x)-T(y)\|\le \|x-y\|$
for all $x,y\in\reals^n$.
Nonexpansive mappings are, by definition, Lipschitz continuous with Lipschitz constant $1$.
$T:\reals^n\rightarrow\reals^n$ is firmly-nonexpansive if
\[
\|T(x)-T(y)\|^2\le \langle x-y,T(x)-T(y)\rangle
\]
for all $x,y\in \mathbb{R}^n$.
Proximal and projection operators are firmly-nonexpansive.

 \subsection{Duality and primal subvalue}
 \label{ss:p-subval}
We call the optimization problem
\begin{equation}
\begin{array}{ll}
\underset{x\in \reals^n}{\mbox{minimize}}&
f(x)+g(x),
\end{array}
\label{eq:P}
\tag{P}
\end{equation}the primal problem.
We call the optimization problem
\begin{equation}
\begin{array}{ll}
\underset{\nu\in \reals^n}{\mbox{maximize}}&
-f^*(\nu)-g^*(-\nu),
\end{array}
\label{eq:D}
\tag{D}
\end{equation}
the dual problem.
Throughout this paper,
we assume $f$ and $g$ are CPC.

\eqref{eq:P} is feasible if $\vz\in \dom f-\dom g$,
strongly infeasible if $\vz\notin \overline{\dom f-\dom g}$,
and weakly infeasible otherwise.
\eqref{eq:P} falls under exactly one of the three cases.
\eqref{eq:P} is infeasible
if it is not feasible.
\eqref{eq:D}
is feasible if $\vz\in \dom (f^*)+\dom (g^*)$,
strongly infeasible if
$\vz\notin \overline{\dom (f^*)+\dom (g^*)}$,
and weakly infeasible otherwise.

%
%
%

We call 
$p^\star=\inf\{f(x)+g(x)\,|\,x\in \reals^n\}$
the primal optimal value and 
$d^\star=\sup\{-f^*(\nu)-g^*(-\nu)\,|\,\nu\in \reals^n\}$
the dual optimal value.
We let $p^\star=\infty$ if \eqref{eq:P} is infeasible and $d^\star=-\infty$ if \eqref{eq:D} is infeasible.
Weak duality, which always holds, states  $d^\star\le p^\star$.
We say strong duality holds between 
 \eqref{eq:P} and \eqref{eq:D},
if $d^\star= p^\star\in[-\infty,\infty]$.
We say \emph{total duality} holds between \eqref{eq:P} and \eqref{eq:D},
if \eqref{eq:P} has a solution, \eqref{eq:D} has a solution, and strong duality holds.


Define the primal subvalue
of \eqref{eq:P} as
\[
p^-=\lim_{\varepsilon\rightarrow 0^+}
\inf_{x,y\in \mathbb{R}^n}\left\{
f(x)+g(y)\,|\,\|x-y\|\le \varepsilon
\right\}.
\] 
The notion of primal subvalue is standard in conic programming \cite{kretschmer1961programmes,yamasaki1969some,luo1997duality,luo2000conic}.
Here, we generalize it to general convex programs.
The following theorem
is well known \cite{rockafellar1974},
although we have not seen it stated exactly in this form.

\begin{theorem}
\label{thm:p-subval}
If $f$ and $g$ are CPC, then $d^\star=p^-\le p^\star$.
\end{theorem}
With Theorem~\ref{thm:p-subval}, 
we can interpret
strong duality as well-posedness of \eqref{eq:P}.
The primal subvalue $p^-$ is the optimal value of 
\eqref{eq:P} achieved with infinitesimal infeasibilities.
When the infinitesimal infeasibilities
provide a non-infinitestimal improvement to the function value,
we can consider  \eqref{eq:P} ill-posed.


\subsection{Douglas--Rachford operator}
Douglas--Rachford splitting (DRS) applied to \eqref{eq:P} is
\begin{align}
x^{k+1/2}&=\Prox_{\gamma f}(z^k)\nonumber\\
x^{k+1}&=\Prox_{\gamma g}(2x^{k+1/2}-z^k)\label{eq:drs}\\
z^{k+1}&=z^k+x^{k+1}-x^{k+1/2}\nonumber
\end{align}
with a  starting point $z^0\in \reals^n$
and a  parameter  $\gamma>0$.
We also express this iteration more concisely as
$z^{k+1}=T_\gamma(z^k)$
where
\begin{equation*}
T_\gamma=
\frac{1}{2}I
+\frac{1}{2}
(2\Prox_{\gamma g}-I)
(2\Prox_{\gamma f}-I).
\label{eq:drs-op}
\end{equation*}
$T_\gamma:\reals^n\rightarrow\reals^n$ is a firmly-nonexpansive operator,
and we interpret DRS as a fixed-point iteration.
Write $T_1$ for $T_\gamma$ with $\gamma=1$.

The standard analysis of DRS assumes total duality, which, again, 
 means \eqref{eq:P} has a solution,  \eqref{eq:D} has a solution, and $d^\star=p^\star$.

\begin{theorem}[Theorem~7.1 and 8.1 of \cite{bauschke2012} and Proposition~4.8 of \cite{eckstein1989}]
Total duality holds between
\eqref{eq:P} and \eqref{eq:D}
if and only if $T_\gamma$ has a fixed point
for some $\gamma>0$.
If 
total duality holds between
\eqref{eq:P} and \eqref{eq:D},
then DRS converges in that
$z^k\rightarrow z^\star$,
where $x^\star=\Prox_{\gamma f}(z^\star)$ is a solution of \eqref{eq:P}.
If total duality does not hold between
\eqref{eq:P} and \eqref{eq:D},
then DRS diverges in that
$\|z^k\|\rightarrow\infty$.
\label{thm:drs-fixed}
\end{theorem}
Theorem~\ref{thm:drs-fixed} is well known,
although the term ``total duality'' is not always used. More often, total duality 
is assumed by instead assuming a saddle point exists for an appropriate Lagrangian.

\subsection{Fixed-point iterations without fixed points}
\label{ss:inf-drs}
Theorem~\ref{thm:drs-fixed} states the DRS iteration has no fixed points under pathologies.
Analyzing fixed-point iterations without fixed points is 
the first part of our pathological analysis.

Let $T:\reals^n\rightarrow\reals^n$ be a firmly-nonexpansive operator.
Write
\[
\ran(I-T)=\{z-T(z)\,|\,z\in \reals^n\}.
\]
Note that
$T$ has a fixed point if and only if 
$\vz\in \ran(I-T)$.
The closure of this set, $\overline{\ran(I-T)}$, is closed and convex \cite{Pazy1971_asymptotic}.
We call 
\[
v=\Pi_{\overline{\ran(I-T)}}(\vz)
\]
the \emph{infimal displacement vector} 
of $T$.
(The term was coined in \cite{bauschke2014generalized}.)
If $T$ has a fixed point, then $v=\vz$, but $v=\vz$ is possible even when $T$ has no fixed point.

The following classical result by Pazy and Baillon et al.\ 
 elegantly characterizes
the asymptotic behavior of fixed-point iterations with respect to $T$.
\begin{theorem}[Theorem~2 of \cite{Pazy1971_asymptotic} and
Corollary~2.3 of \cite{BaillonBruckReich1978_asymptotic}]
\label{thm:pazy}
If $T$ is firmly-nonexpansive and $v$ is its infimal displacement vector,
the iteration $z^{k+1}=T(z^k)$ satisfies
\[
z^k=-kv+o(k),\qquad
z^{k+1}-z^k\rightarrow -v.
\]
\end{theorem}

Theorem~\ref{thm:pazy} is especially powerful when we can concretely characterize $v$.
Recently, Bauschke, Hare, and Moursi published the following elegant formula.
\begin{theorem}[\cite{bauschke2016}]
\label{thm:bauschke}
The infimal displacement vector $v$ of $T_1$, the DRS operator, satisfies
\[
v=\argmin\left\{\|z\|\,|\,z\in\overline{
\dom f-\dom g
}\cap \overline{\dom f^*+\dom g^*}
\right\}.
\]
\end{theorem}
The original result in \cite{bauschke2016} is more general as it applies to the DRS operator of monotone operators.
In Section~\ref{s:infim}, we use Theorem~\ref{thm:bauschke} and the notion of improving directions to
provide a further concrete characterization of $v$.


\section{Theoretical results}
\label{s:theory}
In this section, we present the main theoretical contribution of this paper.
Our analysis requires a generalized notion of improving directions, so we define it first.
Section~\ref{s:infim} analyzes DRS as a fixed-point iteration without fixed points.
Section~\ref{ss:f-val-anal} analyzes DRS as an optimization method that reduces function values.
Section~\ref{ss:shadow}
directly analyzes the evolution
of the $x^{k+1/2}$ and $x^{k+1}$-iterates
of DRS.
Later in Sections~\ref{s:app-drs} and \ref{s:app-admm}
we combine these results 
to analyze the asymptotic behavior of DRS and ADMM applied to pathological convex programs.


While the formula of Theorem~\ref{thm:bauschke} is known,
the use of improving directions 
to concretely characterize the infimal displacement vector $v$ is new.
An improving direction may or may not exist,
and we analyze both cases.
Our analysis shows that 
existence of an improving direction is a key deciding factor in how DRS behaves.

We say $d\in \mathbb{R}^n$ is a \emph{primal improving direction} 
for \eqref{eq:P} if 
\[
\rec f(d)+\rec g(d)<0.
\]
Note 
$\rec f(d)+\rec g(d)=\rec(f+g)(d)$
when \eqref{eq:P} is feasible.
For simplicity, we only consider primal improving directions when \eqref{eq:P} is feasible.
The notion of (primal) improving direction is standard in conic programming \cite{luo1997duality,nesterov1999infeasible,luo2000conic}.
Here, we extend it to general convex programs of the form  \eqref{eq:P}.


If \eqref{eq:P} is feasible 
and there is a primal improving direction, then $p^\star=-\infty$.
To see why, let  $d$ be a  primal improving direction.
Then 
\[
f(x+\alpha d)+g(x+\alpha d)=
\alpha \rec(f+g)(d)+o(\alpha)
\]
for any $x\in \dom f\cap \dom g$
as $\alpha\rightarrow\infty$, and therefore $p^\star=-\infty$.
However, $p^\star=-\infty$ is possible even when \eqref{eq:P} has no improving direction.
We discuss such an example in Section~\ref{s:app-drs}.


Likewise, we say $d'\in \mathbb{R}^n$ is a \emph{dual improving direction} if 
\[
\rec (f^*)(d')+\rec (g^*)(-d')<0.
\]
If \eqref{eq:D} is feasible and there is a
dual improving direction, then $d^\star=\infty$.

%

\subsection{Infimal displacement vector of the DRS operator}
\label{s:infim}
In this section, we provide a further concrete characterization
of the infimal displacement vector $v$.
When
\eqref{eq:P} or \eqref{eq:D} is strongly infeasible,
Theorem~\ref{thm:bauschke} states $v\ne \vz$.
Our contribution is to show $v$ is an improving direction in this case.
For the sake of simplicity, we first analyze $T_1$
and then translate the results to $T_\gamma$ for $\gamma>0$.

We first consider the case where \eqref{eq:P} is feasible
and characterize $v$
based on the primal improving direction
or the absence of it.

\begin{lemma}
\label{lem:key2}
\eqref{eq:P} has an improving direction
if and only if \eqref{eq:D} is strongly infeasible.
Write
\[
d=-\Pi_{\overline{(\dom f^*+\dom g^*})}(\vz).
\]
If \eqref{eq:P} has an improving direction, 
then $d\ne \vz$ and $d$ is an improving direction.
If \eqref{eq:P} has no improving direction, 
then $d=\vz$.
\end{lemma}
\begin{proof}
We first show
\begin{align}
-\Pi_{\overline{(\dom f^*+\dom g^*})}(\vz)
&=
\Prox_{\rec f+\rec g}(\vz).
\label{eq:lem:key}
\end{align}
Let $A$ and $B$ be nonempty convex sets.
The identities of Section~\ref{s:prelim} tell us
\begin{align*}
(\delta_{\overline{A+B}})^*(x)=
\sigma_{\overline{A+B}}(x)=
\sigma_{A+B}(x)
=
\sigma_{A}(x)+\sigma_{B}(x).
\end{align*}
Setting 
$A=\dom f^*$ and $B=\dom g^*$ gives us
\[
 (\delta_{\overline{\dom f^*+\dom g^*}})^*(x)
= \sigma_{\dom f^*}(x)+\sigma_{\dom g^*}(x).
\]
Based on the identities of Section~\ref{s:prelim}, we have
\begin{align*}
\Pi_{\overline{\dom f^*+\dom g^*}}(\vz)
&=
\Prox_{\delta_{\overline{\dom f^*+\dom g^*}}}(\vz)\\
&=
(I-\Prox_{\sigma_{\overline{\dom f^*+\dom g^*}}})(\vz)\\
&=
-\Prox_{\sigma_{\dom f^*+\dom g^*}}(\vz)\\
&=
-\Prox_{\sigma_{\dom f^*}+\sigma_{\dom g^*}}(\vz)\\
&=
-\Prox_{\rec f+\rec g}(\vz).
\end{align*}

Remember that $\rec f+\rec g$ is a convex positively homogeneous function.
Since $\rec f(\vz)+\rec g(\vz)=0$,
\[
\vz=\argmin
\left\{
\rec f(x)+\rec g(x)+(1/2)\|x\|^2
\right\}=\Prox_{\rec f+\rec g}(\vz)
\]
if and only if $\rec f(x)+\rec g(x)\ge 0$ for all $x\in \reals^n$.
By our definition of an improving direction,
 $\rec f(x)+\rec g(x)\ge 0$ for all $x\in \reals^n$
 if and only if there is no improving direction.
By definition, $\vz=\Pi_{\overline{\dom f^*+\dom g^*}}(\vz)$
if and only if \eqref{eq:D} is not strongly infeasible.
So with \eqref{eq:lem:key}, we conclude 
 \eqref{eq:P} has an improving direction if and only if 
\eqref{eq:D} is strongly infeasible.

It remains to show that 
\[
d=\argmin
\left\{
\rec f(x)+\rec g(x)+(1/2)\|x\|^2
\right\}
\]
is an improving direction, if $d\ne \vz$.
Since $d$ is defined as a minimizer, we have
\[
\rec f(d)+\rec g(d)+(1/2)\|d\|^2\le 
\rec f(\vz)+\rec g(\vz)+(1/2)\|\vz\|^2=0.
\]
This implies 
$\rec f(d)+\rec g(d)\le -(1/2)\|d\|^2<0$, i.e., $d$ is an improving direction.
\qed\end{proof}

\begin{lemma}
\label{lem:p-imp-dir}
Assume \eqref{eq:P} is feasible.
Then 
\[
v=-d=\Pi_{\overline{\dom f^*+\dom g^*}}(\vz)
\]
is the infimal displacement vector of $T_{1}$.
\end{lemma}
\begin{proof}
Let $x_0$ be a feasible point of \eqref{eq:P}.
Since $\rec f(d)+\rec g(d)\le 0<\infty$ by Lemma~\ref{lem:key2}
and the definition of an improving direction,
we have $x_0\in \dom  f$, $x_0+d\in \dom g$, and thus $-d\in \dom f-\dom g\subseteq \overline{\dom f-\dom g}$.
Since $-d$ is the minimum-norm element of
$\overline{\dom f^*+\dom g^*}$,
Theorem~\ref{thm:bauschke} tells us that $-d$ is the infimal displacement vector of $T_1$.
\qed\end{proof}

\begin{corollary}
\label{cor:pfea-dinfeas}
Assume \eqref{eq:P} is feasible, and \eqref{eq:D} is feasible. Then $v=\vz$
is the infimal displacement vector of $T_\gamma$
for any $\gamma>0$.
\end{corollary}
\begin{corollary}
\label{cor:pfdwif}
Assume \eqref{eq:P} is feasible, and \eqref{eq:D} is weakly infeasible. Then $v=\vz$
is the infimal displacement vector of $T_\gamma$
for any $\gamma>0$.
\end{corollary}

\begin{corollary}
\label{cor:pfdsif}
Assume \eqref{eq:P} is feasible, and \eqref{eq:D} is strongly infeasible. Then
\[
v=-\gamma d=\gamma \Pi_{\overline{\dom f^*+\dom g^*}}(\vz)\ne \vz
\]
is the infimal displacement vector of $T_\gamma$
for any $\gamma>0$.
Furthermore, $d$ is an improving direction of \eqref{eq:P}.
\end{corollary}

Next, we consider the case where \eqref{eq:D} is feasible
and characterize the infimal displacement vector
based on the dual improving direction
or the absence of it.
\begin{lemma}
\label{lem:pinfeas}
Assume \eqref{eq:D} is feasible. Then 
\[
v=-d'=\Pi_{\overline{\dom f-\dom g}}(\vz)
\]
is the infimal displacement vector of $T_1$.
\end{lemma}

\begin{proof}
Following the same logic as in the proof of Lemma~\ref{lem:key2}, we have
\[
\Pi_{\overline{\dom f-\dom g}}(\vz)
=
-\argmin_\nu\left\{
\rec (f^*)(\nu)+\rec (g^*)(-\nu)+(1/2)\|\nu\|^2
\right\},
\]
and 
\[
d'=-\Pi_{\overline{\dom f-\dom g}}(\vz)
\]
is a dual improving direction, if $d'\ne \vz$.

Let $\nu_0$ be any feasible point of \eqref{eq:D}.
Then 
$\nu_0\in \dom f^*$ and $-\nu_0-d'\in \dom g^*$.
Therefore, $-d'\in \dom f^*+\dom g^*\subseteq \overline{\dom f^*+\dom g^*}$.
Since $-d'$ is defined to be the minimum-norm element of $\overline{\dom f-\dom g}$
we conclude the statement with Theorem~\ref{thm:bauschke}.
\qed\end{proof}

\begin{corollary}
\label{cor:pwinfdf}
Assume \eqref{eq:D} is feasible, and \eqref{eq:P} is weakly infeasible. Then $v=\vz$
is the infimal displacement vector of $T_\gamma$
for any $\gamma>0$.
\end{corollary}

\begin{corollary}
\label{cor:psinfdf}
Assume \eqref{eq:D} is feasible, and \eqref{eq:P} is strongly infeasible. Then
\[
v=-d'=\Pi_{\overline{\dom f-\dom g}}(\vz)\ne \vz,
\]
is the infimal displacement vector of $T_\gamma$
for any $\gamma>0$.
Furthermore, $d'$ is a dual improving direction.
\end{corollary}
Note that for Corollary~\ref{cor:psinfdf}, the infimal displacement vector is independent of the value of $\gamma$.

\subsection{Function-value analysis}
\label{ss:f-val-anal}
In this section, we present the second major theoretical component to our analysis.
Section~\ref{s:infim} analyzed the infimal displacement vector of $T_\gamma$.
This, however, is not sufficient for characterizing the asymptotic
behavior of DRS in relation to the original optimization problem \eqref{eq:P}.

Let us briefly discuss why function-value analysis is necessary.
Consider the convex function
$h(x,y)=x^2/y$
defined for $y>0$.
Note that $h$ has minimizers, $(0,y)$ for any $y>0$,
and the operator $I-\nabla h$ has fixed points.
It is straightforward to verify that $h(\sqrt{y},y)-\inf f\nrightarrow 0$, but $\nabla h(\sqrt{y},y)\rightarrow 0$ as $y\rightarrow\infty$, i.e.,
$(\sqrt{y},y)$ for large $y$ is not an approximate minimizer for $h$
but does approximate satisfy the fixed point condition for $I-\nabla h$.
It is possible to construct a similar example with the DRS operator.
If we let $f=h$ and $g=0$, then DRS reduces to the proximal point method on $h$.
This operator exhibits the same exact issue.

This means \textbf{approximate fixed points do not always correspond to approximate solutions of the original problem}.
This is why we need a separate and distinct function-value analysis to accompany the fixed-point theory.

We now present function-value analysis.
Throughout this section, write
 $x^{k+1/2}$ and $x^{k+1}$ to denote the DRS iterates of \eqref{eq:drs}.

\begin{lemma}
\label{lem:descent-ineq}
For all $k=0,1,\dots$ and any $x\in \reals^n$
\begin{align*}
&f(x^{k+1/2})+g(x^{k+1})
-
f(x)-g(x)\\
&\qquad\qquad\qquad\qquad\le
(1/\gamma)\langle x^{k+1}-x^{k+1/2}, x-z^{k+1}\rangle.
\end{align*}
\end{lemma}
An inequality similar to that of Lemma~\ref{lem:descent-ineq}
has been presented as Proposition~2 of \cite{davis2016}. 
We nevertheless quickly show a direct proof.
\begin{proof}
Write
\begin{align*}
\tilde{\nabla}f(x^{k+1/2})&=(1/\gamma)(z^k-x^{k+1/2})\\
\tilde{\nabla}g(x^{k+1})&=(1/\gamma)(2x^{k+1/2}-z^k-x^{k+1}).
\end{align*}
From the definition of the DRS iteration \eqref{eq:drs},
we can verify that
\begin{align*}
\tilde{\nabla}f(x^{k+1/2})\in\partial f(x^{k+1/2}), \qquad
\tilde{\nabla}g(x^{k+1})\in\partial g(x^{k+1})
\end{align*}
and that
\[
\tilde{\nabla}f(x^{k+1/2})+\tilde{\nabla}g(x^{k+1})=(1/\gamma)(x^{k+1/2}-x^{k+1}).
\]
We also have
\[
z^{k+1}=z^k
-\gamma \tilde{\nabla}f(x^{k+1/2})-\gamma \tilde{\nabla}g(x^{k+1})=x^{k+1/2}-\gamma \tilde{\nabla}g(x^{k+1}).
\]

If $x\notin \dom f\cap \dom g$, then 
$f(x)+g(x)=\infty$ for all $x\in \mathbb{R}^n$
, and there is nothing to prove.
Now, consider any $x\in \dom f\cap \dom g$. Then, by definition of subdifferentials,
\begin{align*}
f(&x^{k+1/2})-f(x)+g(x^{k+1})-g(x)\\
&\le 
\langle \tilde{\nabla}f(x^{k+1/2}), x^{k+1/2}-x\rangle
+
\langle\tilde{\nabla}g(x^{k+1}), x^{k+1}-x\rangle\\
&=
(\tilde{\nabla}f(x^{k+1/2})+\langle\tilde{\nabla}g(x^{k+1}), x^{k+1/2}-x\rangle
+\langle\tilde{\nabla}g(x^{k+1}), x^{k+1}-x^{k+1/2}\rangle\\
&=\langle x^{k+1}-x^{k+1/2}, 
\tilde{\nabla}g(x^{k+1})-(1/\gamma)(x^{k+1/2}-x)\rangle\\
&=(1/\gamma)\langle x^{k+1}-x^{k+1/2}, x-z^{k+1}\rangle.
\end{align*}
\qed\end{proof}

The following result, which is well known for non-pathological setups,
also holds under pathologies, so long as $d^\star=p^\star$.
\begin{lemma}
\label{lem:mean-conv2}
Assume $p^\star=d^\star\in [-\infty,\infty]$.
Assume the infimal displacement vector $v$ of $T_\gamma$ satisfies $v=\vz$.
Then
\[
\lim_{k\rightarrow\infty}
\frac{1}{k+1}\sum^k_{i=0}
f(x^{i+1/2})+g(x^{i+1})=p^\star
\]
and
\[
\liminf_{k\rightarrow\infty}
f(x^{k+1/2})+g(x^{k+1})=
p^\star.
\]
\end{lemma}
\begin{proof}
If $\Delta^0,\Delta^1,\dots$ is any sequence in $\reals^n$,
then
\begin{align*}
\sum^k_{j=0}\langle\Delta^j, \sum^j_{i=0}\Delta^i\rangle &=
\sum^k_{j=0}\sum^k_{i=0}\mathbbm{1}\{i\le j\}
\langle \Delta^j, \Delta^i\rangle\\
&=
\frac{1}{2}\left\|
\sum^k_{i=0}\Delta^i\right\|^2
+
\frac{1}{2}\sum^k_{i=0}\left\|
\Delta^i\right\|^2.
\end{align*}
Let
$\Delta^{k}=z^{k+1}-z^{k}=x^{k+1}-x^{k+1/2}$
and sum the inequality of Lemma~\ref{lem:descent-ineq} to get
\begin{align*}
\gamma
\sum^k_{i=0}
f(x^{i+1/2})&-f(x)+g(x^{i+1})-g(x)
\le
\sum^k_{j=0}\langle\Delta^{j}, x-z^0\rangle
-\sum^k_{j=0}\langle \Delta^{j}, \sum^j_{i=0}\Delta^{i}\rangle\\
&=
\langle z^{k+1}-z^0, x-z^0\rangle-
\frac{1}{2}
\|z^{k+1}-z^0\|^2
-\frac{1}{2}
\sum^k_{i=0}\|z^{i+1}-z^i\|^2\\
&=
-\frac{1}{2}\|z^{k+1}\|^2
+\frac{1}{2}\|z^0\|^2
+\langle z^{k+1}-z^0, x\rangle
-\frac{1}{2}
\sum^k_{i=0}\|z^{i+1}-z^i\|^2.
\end{align*}
Divide both sides by $(k+1)/2$ to get
\begin{align}&
\frac{2\gamma}{k+1}
\sum^k_{i=0}
\left(
f(x^{i+1/2})-f(x)+g(x^{i+1})-g(x) \right)
\label{eq:lem:mean-ineq}
\\
&\quad\le 
-\frac{1}{k+1}\|z^{k+1}\|^2
+\frac{1}{k+1}\|z^{0}\|^2
-\frac{1}{k+1}\sum_{i=0}^k \|z^{i+1}-z^i\|^2
+\frac{2}{k+1}
\langle z^{k+1}-z^0, x\rangle.\nonumber
\end{align}
for all $k=0,1,\dots$ and any $x\in \reals^n$.

We now show
\[
\limsup_{k\rightarrow\infty}
\frac{1}{k+1}
\sum^k_{i=0}
\left(f(x^{i+1/2})+g(x^{i+1})\right)
\le 
p^\star.
\]
Assume $p^\star<\infty$, as otherwise there is nothing to prove.
Let $x$ be any $x\in \dom f\cap \dom g$.
By Theorem~\ref{thm:pazy}, $z^{k}= -k v+o(k)$.
If $v\ne \vz$, then the first (negative) term on the right-hand side of \eqref{eq:lem:mean-ineq} dominates
the positive terms.
If $v= \vz$, then both nonnegative terms 
on the right-hand side of \eqref{eq:lem:mean-ineq} converge to $0$.
In both cases, we have
\begin{equation}
\limsup_{k\rightarrow\infty}
\frac{1}{k+1}
\sum^k_{i=0}
f(x^{i+1/2})+g(x^{i+1})
\le 
f(x)+g(x)
\label{eq:lem:mean-conv}
\end{equation}
for all $x\in \dom f\cap \dom g$.
We minimize the right-hand side to obtain $p^\star$.

By Theorem~\ref{thm:pazy},
$v=\vz$ implies $x^{k+1/2}-x^{k+1}\rightarrow \vz$.
In turn, by Theorem~\ref{thm:p-subval}, 
we have
\[
\liminf_{k\rightarrow\infty}f(x^{k+1/2})+g(x^{k+1})\ge p^\star.
\]
Combining this with \eqref{eq:lem:mean-conv}
gives us the first stated result.

It is straightforward to verify that 
if a real-valued sequence $a^k$ satisfies
\[
\liminf_{k\rightarrow\infty} a^k\ge a,
\qquad
\lim_{k\rightarrow\infty} \frac{1}{k}\sum^k_{i=1}a^i= a,
\]
then
\[
\liminf_{k\rightarrow\infty} a^k= a.
\]
The second stated result follows from this argument.
\qed\end{proof}

Lemma~\ref{lem:mean-conv2} provides the function-value analysis when $v=\vz$,
and the first part of Lemma~\ref{lem:f-val-d} provides the analysis when $v\ne \vz$.
The later parts part of Lemma~\ref{lem:f-val-d} is used
in translating the analyses for DRS to analyses for ADMM
in Section~\ref{s:app-admm}.

\begin{lemma}
\label{lem:f-val-d}
Assume \eqref{eq:P} is feasible and $v\ne \vz$, i.e., \eqref{eq:P} has an improving direction. Then
\[
f(x^{k+1/2})+g(x^{k+1})\rightarrow p^\star=-\infty.
\]
Moreover, $|f(x^{k+1/2})|\le\mathcal{O}(k)$ and $|g(x^{k+1})|\le \mathcal{O}(k)$
as $k\rightarrow \infty$.
Assume \eqref{eq:P} is feasible and $v=\vz$.
Then $|f(x^{k+1/2})|\le o(k)$ and $|g(x^{k+1})|\le o(k)$
as $k\rightarrow \infty$.
\end{lemma}
\begin{proof}
When \eqref{eq:P} has an improving direction,
Corollary~\ref{cor:pfdsif} and Theorem~\ref{thm:pazy} tells us 
\[
z^{k+1}-z^k=x^{k+1}-x^{k+1/2}\rightarrow \gamma d.
\]
Then Lemma~\ref{lem:descent-ineq} tells us that
\begin{align}
(1/k)(f(x^{k+1/2})+g(x^{k+1}))
&\le 
-
(1/\gamma)\langle x^{k+1}-x^{k+1/2}, (1/k)z^{k+1}\rangle
+O(1/k)\nonumber
\end{align}
which tells us
\begin{equation}
\limsup_{k\rightarrow\infty}(1/k)(f(x^{k+1/2})+g(x^{k+1}))
\le 
-
\gamma \|d\|^2.
\label{eq:fill}
\end{equation}
This proves the first statement.

Assume $v=\vz$.
With the same reasoning as for \eqref{eq:fill} we get
\[
\limsup_{k\rightarrow\infty}(1/k)(f(x^{k+1/2})+g(x^{k+1}))\le 0.
\]

Assume  \eqref{eq:P} feasible, without making any asumptions on $v$.
Write $\tilde{\nabla} f(x^{1/2})$ for any subgradient of $f$ at $x^{1/2}$. Then
\begin{align*}
f(x^{k+1/2})&\ge 
f(x^{1/2})+
\langle \tilde{\nabla} f(x^{1/2}), x^{k+1/2}-x^{1/2}\rangle\\
&\ge f(x^{1/2})-
\|\tilde{\nabla} f(x^{1/2})\|
\|x^{k+1/2}-x^{1/2}\|
=
k\gamma \|d\|\|\tilde{\nabla} f(x^{1/2})\|
+o(k),
\end{align*}
and we conclude
\[
\liminf_{k \rightarrow \infty}(1/k)f(x^{k+1/2})\ge -\gamma \|d\|\|\tilde{\nabla} f(x^{1/2})\|.
\]
With a similar argument, we get
\[
\liminf_{k \rightarrow \infty}(1/k)g(x^{k+1})\ge -\gamma \|d\|\|\tilde{\nabla} g(x^{1})\|
\]
where $\tilde{\nabla} g(x^{1})$ is any subgradient of $g$ at $x^1$.
Combining these with \eqref{eq:fill} gives us the remaining statements.
%
%
%
%
%
%
%
\qed\end{proof}

\begin{lemma}
\label{lem:if-conv}
Assume $p^\star=d^\star$.
Assume $x^{k+1/2}$ and $x^{k+1}$ are the DRS iterates 
as defined in \eqref{eq:drs}.
If $x^{k+1/2},x^{k+1}\rightarrow x^\star$ for some $x^\star$, then $x^\star$ is a solution.
\end{lemma}
\begin{proof}
We first note that closed functions are by definition lower semi-continuous,
and that $f$ and $g$ are assumed to be closed.
By Lemma~\ref{lem:mean-conv2} we have
\begin{align*}
f(x^\star)+g(x^\star)&\le 
\liminf_{k\rightarrow \infty}
f(x^{k+1/2})+g(x^{k+1})=p^\star,
\end{align*}
and we conclude $f(x^\star)+g(x^\star)=p^\star$.
\qed\end{proof}

\subsection{Evolution of shadow iterates}
\label{ss:shadow}
Section~\ref{s:infim}
characterized the evolution of the $z^k$-iterates,
which we could call the main iterates.
The $x^{k+1/2}$ and $x^{k+1}$-iterates of DRS
are  called
the \emph{shadow iterates}.
Here, we analyze the evolution
of the shadow iterates.

Although the results of this section are
are not as fundamental or important
as the results of Sections~\ref{s:infim} and \ref{ss:f-val-anal},
we do need these results later, especially when translating the analyses for DRS to analyses for ADMM.


\begin{lemma}
\label{lem:shadow}
If $v=\vz$, then 
 $x^{k+3/2}-x^{k+1/2}\rightarrow \vz$
and
$x^{k+2}-x^{k+1}\rightarrow \vz$.
\end{lemma}
\begin{proof}
Since $v=\vz$, we have $z^{k+1}-z^k\rightarrow \vz$.
Since the map the defines $z^{k}\mapsto x^{k+1/2}$ 
and $z^{k+1}\mapsto x^{k+3/2}$ 
is Lipschitz continuous,
$x^{k+3/2}-x^{k+1/2}\rightarrow \vz$.
Finally,
$z^{k+1}-z^k\rightarrow \vz$ and 
$x^{k+3/2}-x^{k+1/2}\rightarrow \vz$
implies
$x^{k+2}-x^{k+1}\rightarrow \vz$.
\qed\end{proof}

\begin{lemma}
\label{lem:shadow-pinfeas}
If (P) is strongly infeasible and (D) is feasible,
then  $x^{k+3/2}-x^{k+1/2}\rightarrow \vz$
and
$x^{k+2}-x^{k+1}\rightarrow \vz$.
\end{lemma}
\begin{proof}
Write $-d'$ for the infimal displacement vector
as given by Corollary~\ref{cor:psinfdf}.
By Theorem~\ref{thm:pazy}, we have
\[
z^{k+1}-z^k=x^{k+1}-x^{k+1/2}\rightarrow d'.
\]
The projection inequality states 
\begin{equation}
\langle v-\Pi_C x,  \Pi_C x -x\rangle \ge 0
\label{eq:proj}
\end{equation}
for any nonempty closed convex set $C$, $v\in C$, and $x\in \mathbb{R}^n$.
Since $-d'=\Pi_{\overline{\dom f-\dom g}}(\vz)$,
\eqref{eq:proj} tells us that
\[
\langle
d',x-x^{k+1}\rangle +\|d'\|^2\le 0
\]
for any $x\in \dom f$.
Using  
$x^{k+1/2}=\Prox_{\gamma f}(z^k)$
and firm-nonexpansiveness of $\Prox$, we get
\begin{align*}
\|\Prox_{\gamma f}(z^k+d')-x^{k+1/2}\|^2
&\le 
\langle d',\Prox_{\gamma f}(z^k+d')-x^{k+1/2}\rangle\\
&=
\langle d',\Prox_{\gamma f}(z^k+d')-x^{k+1}\rangle
+\langle d',x^{k+1}-x^{k+1/2}\rangle\\
&\rightarrow 0
\end{align*}
since $\langle d',x^{k+1}-x^{k+1/2}\rangle\rightarrow \|d'\|^2$.
So $\Prox_{\gamma f}(z^k+d')-\Prox_{\gamma f}(z^k)\rightarrow \vz$.
Since $\Prox$ is Lipschitz continuous, 
$z^{k+1}-z^k-d'\rightarrow \vz$ implies
\[
\Prox_{\gamma f}(z^{k}+d')-
\Prox_{\gamma f}(z^{k+1})\rightarrow \vz.
\]
Putting everything together we conclude
\[
\Prox_{\gamma f}(z^{k+1})-
\Prox_{\gamma f}(z^k)
=
x^{k+3/2}-x^{k+1/2}
\rightarrow \vz.
\]

Since 
\[
z^{k+2}-z^{k+1}
=\underbrace{z^{k+1}-z^k}_{\rightarrow d'}
+x^{k+2}-x^{k+1}-
\underbrace{(x^{k+3/2}-x^{k+1/2})}_{\rightarrow\vz}
\rightarrow d'
\]
we also conclude that $x^{k+2}-x^{k+1}\rightarrow \vz $.
\qed\end{proof}

\begin{lemma}
\label{lem:shadow3}
If (P) has an improving direction, 
and (P) is feasible, then $x^{k+3/2}-x^{k+1/2}\rightarrow \gamma d$
and
$x^{k+2}-x^{k+1}\rightarrow \gamma d$,
where $-\gamma d=\gamma \Pi_{\overline{\dom f^*+\dom g^*}}(\vz)$ 
is the infimal displacement vector as given in Corollary~\ref{cor:pfdsif}.
\end{lemma}

\begin{proof}
For simplicity, assume $\gamma=1$.
For $\gamma \ne 1$, we scale $f$ and $g$ to get the stated result.

Rewrite the DRS iteration as
\begin{align*}
x^{k+1/2}&=\Prox_{f}(z^k)\\
\nu^{k+1/2}&=z^k-x^{k+1/2}=\Prox_{f^*}(z^k)\\
x^{k+1}&=\Prox_{ g}(2x^{k+1/2}-z^k)\\
\nu^{k+1}&=2x^{k+1/2}-z^k-x^{k+1}=\Prox_{g^*}(2x^{k+1/2}-z^k)\\
z^{k+1}&=z^k-(\nu^{k+1}+\nu^{k+1/2}).
\end{align*}
By Theorem~\ref{thm:pazy}, we have
\[
z^{k+1}-z^k=x^{k+1}-x^{k+1/2}\rightarrow d.
\]
By the same reasoning as in Lemma~\ref{lem:shadow-pinfeas},
we can use \eqref{eq:proj} and firm-nonexpansiveness
to show that 
\[
\nu^{k+3/2}-\nu^{k+1/2}=
\Prox_{f^*}(z^{k+1})-\Prox_{f^*}(z^k)\rightarrow \vz.
\]
Since 
\[
z^{k+1}-z^k=
\underbrace{\nu^{k+3/2}-\nu^{k+1/2}}_{\rightarrow \vz}+x^{k+3/2}-x^{k+1/2}\rightarrow d,
\]
we have
$x^{k+3/2}-x^{k+1/2}\rightarrow d$.

Since 
\[
z^{k+2}-z^{k+1}
=\underbrace{z^{k+1}-z^k}_{\rightarrow d}
+x^{k+2}-x^{k+1}-
\underbrace{(x^{k+3/2}-x^{k+1/2})}_{\rightarrow d}
\rightarrow d
\]
we also conclude that $x^{k+2}-x^{k+1}\rightarrow d $.
\qed\end{proof}

%
%

\section{Pathological convergence: DRS}
\label{s:app-drs}
In this section, we use the theory of Section~\ref{s:theory} to analyze DRS under pathologies.
We classify the status of 
 \eqref{eq:P} and \eqref{eq:D} into 7 cases
and provide convergence analyses for the first 6 cases, the ones that assume strong duality.

\subsection{Classification}
\label{ss:classification}
The primal-dual problem pair, \eqref{eq:P} and \eqref{eq:D},
falls under exactly one of the following 7 distinct cases.
\paragraph{Case (a)}
Total duality holds between \eqref{eq:P} and \eqref{eq:D}.

In other words, \eqref{eq:P} and \eqref{eq:D} have solutions, and $d^\star=p^\star$.
For example, the primal problem
\begin{equation*}
\begin{array}{ll}
\mbox{minimize}& x-\log x
\end{array}
\end{equation*}
and its dual problem
\[
\begin{array}{ll}
\mbox{maximize}& 1+\log(y)\\
\mbox{subject to} &y=1
\end{array}
\]
both have solutions, and $d^\star=p^\star=1$.

\paragraph{Case (b)}
$d^\star=p^{\star}$ is finite, \eqref{eq:P} has a solution, \eqref{eq:D} has no solution.

For example, the primal problem

\[
\begin{array}{ll}
\mbox {minimize} & \underbrace{\delta_{\{(x_1,x_2)\,|\,x_1^2+x_2^2\leq 1\}}(x_1,x_2)}_{f(x)}+\underbrace{x_2+\delta_{\{(x_1,x_2)\,|\,x_1=1\}}(x_1,x_2)}_{g(x)}
\end{array}
\]
has a solution but its dual problem
\begin{equation*}
\begin{array}{ll}
\mbox{maximize}& -\sqrt{\nu_1^2+\nu_2^2}+\nu_1-\delta_{\{\nu_2=1\}}(-\nu_2)
\end{array}
\end{equation*}
does not. Nevertheless, $d^\star=p^\star=0$.
%

\paragraph{Case (c)}
$d^\star=p^{\star}$ is finite, \eqref{eq:P} is feasible, 
but \eqref{eq:P} has no solution. 

To get such an example,
swap the role of the primal and the dual in the example for case (b).

\paragraph{Case (d)}
$d^\star=p^{\star}=-\infty$,
\eqref{eq:P} is feasible, 
but there is no improving direction.

This implies \eqref{eq:D} is weakly infeasible.
For example, the primal problem
\begin{equation*}
\begin{array}{ll}
\mbox{minimize}& \delta_{\{x\,|\,x\ge 1\}}(x)-\log x
\end{array}
\end{equation*}
has no solution and has optimal value $p^\star=-\infty$.
Since the derivative of the objective, $-1/x$, goes to $0$ as $x\rightarrow \infty$,
the primal problem has no improving direction.
The dual problem
\[
\begin{array}{ll}
\mbox{maximize}& y+1+\log(y)\\
\mbox{subject to} &y\le 0
\end{array}
\]
is weakly infeasible.

\paragraph{Case (e)}
$d^\star=p^{\star}=-\infty$, \eqref{eq:P} is feasible, 
and there is an improving direction.

This implies \eqref{eq:D} is strongly infeasible.
For example, the primal problem
\[
\begin{array}{ll} 
\mbox{minimize}& x+x
\end{array}
\]
has an improving direction, namely $d=-1$,
and the dual problem
\[
\begin{array}{ll}
\mbox{maximize}& \delta_{\{1\}}(x)+\delta_{\{1\}}(-x)
\end{array}
\]
is strongly infeasible.


\paragraph{Case (f)}
$d^\star=p^{\star}=\infty$ and \eqref{eq:P} is infeasible.

For example, the problem
\begin{equation*}
\begin{array}{ll}
\mbox{minimize}& 
1/\sqrt{-x}-\log(x)
\end{array}
\end{equation*}
is infeasible, and its dual
\begin{equation*}
\begin{array}{ll}
\mbox{maximize}& 
(3/2^{2/3})y^{1/3}+1+\log(y)
\\
\mbox{subject to}&
y\ge 0
\end{array}
\end{equation*}
has optimal value $d^\star=\infty$.

\paragraph{Case (g)}
$d^\star<p^{\star}$, i.e. strong duality fails.

\subsection{Convergence results}
\label{ss:drs-conv}

\begin{theorem}\cite{lions1979,davis2016}
\label{thm:drs-a}
In case (a),
$x^{k+1/2},x^{k+1}\rightarrow x^\star$, where $x^\star$ is a solution of \eqref{eq:P}
and
\[
\lim_{k\rightarrow\infty}
f(x^{k+1/2})+g(x^{k+1})=p^\star.
\]
\end{theorem}

\begin{theorem}
\label{thm:drs-b}
In case (b), $x^{k+1}-x^{k+1/2}\rightarrow \vz$ and
\[
\lim_{k\rightarrow\infty}
\frac{1}{k+1}\sum^k_{i=0}
f(x^{i+1/2})+g(x^{i+1})=p^\star,
\qquad
\liminf_{k\rightarrow\infty}f(x^{k+1/2})+g(x^{k+1})=p^\star.
\]
Furthermore, if $x^{k+1/2}\rightarrow x^\star$
(or equivalently if $x^{k+1}\rightarrow x^\star$)
then $x^\star$ is a solution.
\end{theorem}
\begin{proof}
This follows from Theorem~\ref{thm:pazy}, Corollary~\ref{cor:pfea-dinfeas}, Lemma~\ref{lem:mean-conv2}, and Lemma~\ref{lem:if-conv}.
\qed\end{proof}

\begin{theorem}
\label{thm:drs-c}
In case (c), $x^{k+1}-x^{k+1/2}\rightarrow \vz$,
\[
\lim_{k\rightarrow\infty}
\frac{1}{k+1}\sum^k_{i=0}
f(x^{i+1/2})+g(x^{i+1})=p^\star,
\qquad
\liminf_{k\rightarrow\infty}f(x^{k+1/2})+g(x^{k+1})=p^\star,
\]
and $(x^{k+1/2},x^{k+1})$ do not converge.
\end{theorem}
\begin{proof}
This follows from Theorem~\ref{thm:pazy}, Corollary~\ref{cor:pfea-dinfeas}, Lemma~\ref{lem:mean-conv2}, and the contrapositive of Lemma~\ref{lem:if-conv}.
\qed\end{proof}

\begin{theorem}
\label{thm:drs-d}
In case (d),
\eqref{eq:D} is weakly infeasible,
$x^{k+1}-x^{k+1/2}\rightarrow \vz$,
\[
\lim_{k\rightarrow\infty}
\frac{1}{k+1}\sum^k_{i=0}
f(x^{i+1/2})+g(x^{i+1})=-\infty,
\qquad
\liminf_{k\rightarrow\infty}f(x^{k+1/2})+g(x^{k+1})=-\infty,
\]
and $(x^{k+1/2},x^{k+1})$ do not converge.
\end{theorem}
\begin{proof}
This follows from Theorem~\ref{thm:pazy}, Lemma~\ref{lem:key2},  Corollary~\ref{cor:pfdwif}, Lemma~\ref{lem:mean-conv2}, and the contrapositive of Lemma~\ref{lem:if-conv}.
\qed\end{proof}

\begin{theorem}
\label{thm:drs-e}
In case (e), 
\eqref{eq:D} is strongly infeasible,
$x^{k+1}-x^{k+1/2}\rightarrow \gamma d$,
where $d$ is an improving direction,
\[
\lim_{k\rightarrow\infty}f(x^{k+1/2})+g(x^{k+1})=-\infty,
\]
and $(x^{k+1/2},x^{k+1})$ do not converge.
Furthermore,
$\dist(x^{k+1/2},\dom g)\rightarrow 0$
and\\
$\dist(x^{k+1},\dom f)\rightarrow 0$.
\end{theorem}
\begin{proof}
All but the last assertions 
 follows from Theorem~\ref{thm:pazy}, Lemma~\ref{lem:key2}, 
Corollary~\ref{cor:pfdsif}, Lemma~\ref{lem:f-val-d},
and the contrapositive of Lemma~\ref{lem:if-conv}.
By Lemma~\ref{lem:shadow3}
$x^{k+1/2}-x^{k-1/2}\rightarrow \gamma d$
and by 
Theorem~\ref{thm:pazy} and
Corollary~\ref{cor:pfdsif}
$x^{k}-x^{k-1/2}\rightarrow \gamma d$.
So
$x^{k+1/2}-x^k\rightarrow \vz$.
Since $x^{k}\in \dom g$, we have
\[
\dist(x^{k+1/2},\dom g)\le
\dist(x^{k+1/2},x^k)\rightarrow 0.
\]
Since $x^{k+1/2}\in \dom f$, we have
\[
\dist(x^{k},\dom f)\le
\dist(x^{k},x^{k+1/2})\rightarrow 0.
\]
\qed\end{proof}
\begin{theorem}
\label{thm:drs-f}
In case (f), $\|x^{k+1}-x^{k+1/2}\|\rightarrow \dist(\dom f,\dom g)$.
\end{theorem}
\begin{proof}
This follows from Theorem~\ref{thm:pazy} and Corollaries~\ref{cor:pwinfdf} and \ref{cor:psinfdf}.
\qed\end{proof}

\subsection{Interpretation}
\label{ss:interpretation}
We can view the DRS 
as an algorithm with two major goals: make the iterates feasible and optimal.
With some caveats, DRS succeeds at both.
As an auxiliary goal, we want the shadow iterates of DRS
to converge to a solution if one exists.
With some caveats, DRS succeeds at this as well.
Finally, DRS provides a certificate of infeasibility in cases (e) and (f).

In cases (a), (b), (c), and (d)
the iterates become approximately feasible in that $x^{k+1}-x^{k+1/2}\rightarrow \vz$.
In case (e)
the iterates become approximately feasible in that $\dist(x^{k+1/2},\dom g)\rightarrow 0$ and
$\dist(x^{k+1},\dom f)\rightarrow 0$.
In case (f), feasibility is impossible,
but DRS does its best to achieve feasibility.

In cases (a), (b), (c), (d), and (e),
the function values on average converge to the optimal value.
In other words, DRS finds the correct optimal value in these cases.

In case (a), the shadow iterates, the $x^{k+1/2}$ and $x^{k+1}$ iterates, converge to a solution.
In case (b), we do not know whether the shadow iterates converge to a solution.
However, if they converge, the limit is a solution.
In cases (c), (d), and (e), the shadow iterates do not converge, which is good
since there is no solution to converge to.

In cases (e) and (f), the limit $z^{k+1}-z^k\rightarrow -v\ne \vz$ provides a certificate of dual and primal strong infeasibility, respectively.
These may be computationally useful when verifying the validity of a certificate is easy,
which is the case for conic programs.

We quickly clarify the contribution.
The analysis of case (a) is well known and is not the focus of this work, but we include it's discussion here for completeness.
Approximate feasibility in cases (a), (b), (c) and (d) directly follows from prior work, in particular from Theorems~\ref{thm:pazy} and \ref{thm:bauschke}.
The approximate feasibility results for cases (e) and (f) are contributions of this work.

\subsection{Feasibility problems}
\label{ss:feas}
Consider the problem of finding an $x\in A\cap B$, where $A$ and $B$ are nonempty closed convex sets.
Recasting this convex feasibility problem into an equivalent optimization problem and using
Theorem~\ref{thm:bauschke} \cite{bauschke2016}, 
Theorem~\ref{thm:pazy} \cite{Pazy1971_asymptotic,BaillonBruckReich1978_asymptotic},
Theorem~\ref{thm:drs-a} \cite{lions1979},
and basic convex analysis provides us the following results:
\begin{itemize}
    \item Case (a).
If $A \cap B \ne \emptyset$ 
then $x^{k+1/2},x^{k+1}\rightarrow x^\star$ where $x^\star\in A\cap B$.
\item Case (f).
If $\dist(A, B)>0$, then $\|x^{k+1}-x^{k}\|\rightarrow\dist(A,B)$.
\item Case (g).
If $A \cap B \ne \emptyset$ but $\dist(A, B)=0$, then $x^{k+1/2}-x^{k+1}\rightarrow 0$.
\end{itemize}

Specifically, one can recast the convex feasibility problem  $x\in A\cap B$ into the primal problem
\[
\begin{array}{ll}
\underset{x\in \reals^n}{\mbox{minimize}}&
\delta_A(x)+\delta_B(x),
\end{array}
\]
which has the dual problem
\[
\begin{array}{ll}
\underset{\nu\in \reals^n}{\mbox{maximize}}&
-\sigma_A(\nu)-\sigma_B(-\nu).
\end{array}
\]
When $A\cap B\ne\emptyset$, then $p^\star=0$ with $x\in A\cap B$ and $d^\star=0$ with $\nu=0$.
Therefore total duality holds (i.e., we have case (a)) and Theorem~\ref{thm:drs-a} applies.
When $\dist(A, B)>0$, then $p^\star=\infty$ since $A\cap B=\emptyset$.
For the dual, define $\tilde{\nu}=P_{\overline{A-B}}(0)$, which satisfies 
\[
\langle
a-b
,\tilde{\nu}
\rangle\ge \|\tilde{\nu}\|^2
\]
for all $a\in A$ and $b\in B$ by the optimality conditions defining the projection.
Then we have
\[
-\sigma_A(-\eta\tilde{\nu})-\sigma_B(+\eta\tilde{\nu})
=
\inf_{a\in A,b\in B}\langle a-b,\tilde{\nu}\rangle
\ge 
\eta\|\tilde{\nu}\|^2
\]
for $\eta>0$. 
Since $\|\tilde{\nu}\|=\dist(A, B)>0$, with $\eta\rightarrow\infty$ we conclude $d^\star=\infty$.
So we have case (f) and Theorem~\ref{thm:drs-f} applies.
However, the results of this work say nothing for case (g).
The contribution of this work is to consider improving directions and function-value analysis, but both notions are not relevant in the setup of convex feasibility problems.
Therefore, our work does not provide any new results for the convex feasibility problems.








Prior work on the convex feasibility setup provides further stronger results.
By \cite[Theorem 3.13]{BauschkeCombettesLuke2004_finding}, we have \[
x^{k+1}-x^{k+1/2}\rightarrow \Pi_{\overline{B-A}}(\vz).
\]
Furthermore, by \cite[Theorem 4.5]{BauschkeMoursi2017_douglasrachford}, we have
\[
(x^{k+1/2},x^{k+1})\rightarrow (a^\mathrm{apx},b^\mathrm{apx})\in
\argmin_{(a,b)\in A\times B}\{\|a-b\|\}
\]
if the $\argmin $ is nonempty.
(The pairs in the $\argmin$ are called ``best approximation pairs'' between $A$ and $B$.)
These results show that the relevant dichotomy is whether a best approximation pair exists, rather than whether strong duality holds.
These results cannot be obtained from the analysis of our work.

\section{Pathological convergence: ADMM}
\label{s:app-admm}

We now  analyze ADMM under pathologies.
Consider the primal problem
\begin{equation}
\begin{array}{ll}
\underset{x\in \reals^p,y\in \reals^q}{\mbox{minimize}}& f(x)+g(y)\\
\mbox{subject to}& Ax+By=c,
\end{array}
\tag{P-ADMM}
\label{eq:padmm}
\end{equation}
where $f:\reals^p\rightarrow \reals\cup\{\infty\}$ and
$g:\reals^q\rightarrow \reals\cup\{\infty\}$ are CPC,
$A\in \reals^{n\times p}$, $B\in \reals^{n\times q}$, and $c\in \reals^n$,
and its dual problem
\begin{equation}
\begin{array}{ll}
\underset{\nu\in \reals^n}{\mbox{maximize}}& -f^*(-A^T\nu)-g^*(-B^T\nu)-c^T\nu.
\end{array}
\tag{D-ADMM}
\label{eq:dadmm}
\end{equation}
Write $p^\star$ and $d^\star$ for the primal and dual optimal values.
ADMM applied to this primal-dual problem pair is
\begin{align}
x^{k+1}&\in\argmin_{x\in \reals^p}\left\{f(x)+
\langle \nu^k, Ax+By^k-c\rangle
+
\frac{1}{2\gamma}\|Ax+By^{k}-c\|^2\right\}\nonumber\\
y^{k+1}&\in\argmin_{y\in\reals^q} \left\{g(y)+
\langle \nu^k, Ax^{k+1}+By-c\rangle+
\frac{1}{2\gamma}\|Ax^{k+1}+By-c\|^2\right\}\nonumber\\
\nu^{k+1}&=\nu^k+(1/\gamma)(Ax^{k+1}+By^{k+1}-c)\label{eq:admm}.
\end{align}

For ADMM to be well-defined, 
the argmins of \eqref{eq:admm} must exist.
Throughout this section, we furthermore assume the regularity conditions
\begin{align}
(\ran A^T)&\cap \ri\dom(f^*)\ne \emptyset,\label{eq:rcf}\\
(\ran B^T)&\cap \ri\dom(g^*)\ne \emptyset.\label{eq:rcg}
\end{align}
Here, $\ri$ denotes the relative interior of a set.
These conditions ensure the subproblems are solvable \cite[Theorem~16.3]{rockafellar1970}.

Without these regularity conditions,
the subproblems of  \eqref{eq:admm} may not have solutions.
This is often overlooked and sometimes even misunderstood
throughout the ADMM literature.
(The highly influential paper \cite{boyd2011}
mistakenly claimed it is enough for $f$ and $g$ to be CPC.
Chen, Sun, and Toh  \cite{Chen2017} pointed out that additional assumptions are needed.)

\subsection{Classification and convergence results}
\label{ss:class-conv}
Under \eqref{eq:rcf} and \eqref{eq:rcg},
the status of \eqref{eq:padmm} and \eqref{eq:dadmm}
falls under exactly one of the following 5 distinct cases.

\paragraph{Case (a)}
$d^\star=p^{\star}$, both \eqref{eq:padmm} and \eqref{eq:dadmm} have solutions.
\begin{theorem}[\cite{glowinski1975,boyd2011,davis2016}]
In case (a), $Ax^{k}+By^{k}-c\rightarrow \vz$ and 
\[
\lim_{k\rightarrow\infty}f(x^{k})+g(y^k)\rightarrow p^\star.
\]
\end{theorem}

\paragraph{Case (b)}
$d^\star=p^{\star}$, \eqref{eq:padmm} has a solution,  \eqref{eq:dadmm} has no solution.

\begin{theorem}
\label{thm:admm1}
In case (b), $Ax^{k}+By^{k}-c\rightarrow \vz$ and
\[
\lim_{k\rightarrow\infty}\frac{1}{k}\sum^k_{i=1}f(x^{i})+g(y^i)=p^\star,
\quad
\liminf_{k\rightarrow\infty}f(x^k)+g(y^k)=p^\star.
\]
Furthermore, if $(x^k,y^k)\rightarrow (x^\star, y^\star)$, then $(x^\star, y^\star)$ is a solution.
\end{theorem}

\paragraph{Case (c)}
$d^\star=p^{\star}\in [-\infty,\infty)$, \eqref{eq:padmm} is feasible but has no solution. 

\begin{theorem}
\label{thm:admm2}
In case (c), $Ax^{k}+By^{k}-c\rightarrow \vz$ and
\[
\lim_{k\rightarrow\infty}\frac{1}{k}\sum^k_{i=1}f(x^{i})+g(y^i)=p^\star,
\quad
\liminf_{k\rightarrow\infty}f(x^k)+g(y^k)=p^\star,
\]
and the sequence $(x^k,y^k)$ does not converge.
\end{theorem}





\paragraph{Case (d)}
$d^\star=p^{\star}=\infty$, \eqref{eq:padmm} is infeasible.
\begin{theorem}
\label{thm:admm5}
In case (d), 
\[
\|Ax^{k}+By^{k}-c\|\rightarrow 
\inf_{\substack{x\in \dom f\\y\in \dom g}}\|Ax+By-c\|.
\]
\end{theorem}

\paragraph{Case (e)}
$d^\star<p^{\star}$, i.e. strong duality fails.

\subsection{Interpretation}
\label{ss:admm-interpretation}
With some caveats, ADMM
succeeds at achieving feasibility and optimality.
In cases (a), (b), and (c)
the iterates become approximately feasible in that $Ax^{k}+By^{k}-c\rightarrow \vz$,
and the function values on average converge to the solution.
In case (d), feasibility is impossible, but ADMM does its best to achieve feasibility.

%
%

\subsection{Proofs}
\label{ss:admm-deriv}
ADMM is often analyzed as DRS applied to \eqref{eq:dadmm} \cite{gabay1983}.
In this proof, however, we take the less common approach shown in 
\cite{eckstein1989,yan2016},
which derives ADMM directly from the primal problem.
We do so as the function-value analysis 
of Section~\ref{ss:f-val-anal} translate nicely
with this primal approach.




Consider the equivalent primal optimization problem
\begin{equation*}
\begin{array}{ll}
\underset{z\in \reals^n}{\mbox{minimize}}&
\tilde{f}(z)+\tilde{g}(z)
\end{array}
\end{equation*}
with
\[
\tilde{f}(z)=\inf \{f(x)\,|\,Ax+z=\vz\},\quad
\tilde{g}(z)=\inf \{g(y)\,|\,By-c=z\},
\]
which are CPC functions, as we assume
\eqref{eq:rcf} and \eqref{eq:rcg}
\cite[Theorem~16.3]{rockafellar1970}.
We apply DRS to this form to get
\begin{align*}
\tilde{x}^{k+1/2}&=\argmin_{\tilde{x}}\left\{\gamma\tilde{g}(\tilde{x})+(1/2)\|\tilde{x}-z^k\|^2\right\}\\
\tilde{x}^{k+1}&=\argmin_{\tilde{x}}\left\{\gamma\tilde{f}(\tilde{x})+(1/2)\|\tilde{x}-2\tilde{x}^{k+1/2}+z^k\|^2\right\}\\
z^{k+1}&=z^k+\tilde{x}^{k+1}-\tilde{x}^{k+1/2},
\end{align*}
where we perform the $\tilde{g}$-update before the $\tilde{f}$-update.
We introduce and substitute the variables $x^k$, $y^k$, and $\nu^k$ defined implicitly by
$\tilde{x}^{k+1/2}=By^{k+1}-c$, $\tilde{x}^{k+1}=-Ax^{k+2}$, and
$z^k=-\gamma \nu^k-Ax^{k+1}$ to get
\begin{align*}
y^{k+1}&=\argmin_{y}\left\{\gamma g(y)
+\gamma\langle \nu^k, Ax^{k+1}+By-c\rangle+
(1/2)\|Ax^{k+1}+By-c\|^2\right\}\\
x^{k+2}&=\argmin_{x}\left\{\gamma f(x)+
\gamma \langle \nu^{k+1}, Ax+By^{k+1}-c\rangle +
(1/2)\|Ax+By^{k+1}-c
\|^2\right\}\\
 \nu^{k+1}&= \nu^k+(1/\gamma)(Ax^{k+1}+By^{k+1}-c).
\end{align*}
Reordering the updates to get the dependency right, we get
\begin{align*}
y^{k+1}&=\argmin_{y}\left\{\gamma g(y)
+\gamma\langle \nu^k, Ax^{k+1}+By-c\rangle+
(1/2)\|Ax^{k+1}+By-c\|^2\right\}\\
 \nu^{k+1}&= \nu^k+(1/\gamma)(Ax^{k+1}+By^{k+1}-c)\\
x^{k+2}&=\argmin_{x}\left\{\gamma f(x)+\gamma \langle \nu^{k+1}, Ax+By^{k+1}-c\rangle+(1/2)\|Ax+By^{k+1}-c\|^2\right\}.
\end{align*}
Finally, 
redefine the start and end of an iteration
so that it updates $x^{k+1}$, $y^{k+1}$, and $\nu^{k+1}$
instead $y^{k+1}$, $\nu^{k+1}$, and $x^{k+2}$.
With this, we get \eqref{eq:admm}.



The the last step, where we redefine the start and end of an iteration,
introduces a subtlety
when translating the results of Section~\ref{ss:drs-conv}.
In particular, the results of Section~\ref{ss:shadow} are necessary because of this.


Theorem~\ref{thm:admm1}
follows from Theorem~\ref{thm:drs-b}
and Lemmas~\ref{lem:p-imp-dir},
\ref{lem:f-val-d},
 and \ref{lem:shadow}.
Theorem~\ref{thm:admm5}
follows from 
Theorem~\ref{thm:drs-f} and Lemma~\ref{lem:shadow-pinfeas}.

Case (c) of this section
corresponds to cases (c), (d), and (e)
of Section~\ref{ss:drs-conv}.
For the three cases,
we use
Theorem~\ref{thm:drs-c} and 
Lemmas~\ref{lem:p-imp-dir}, 
\ref{lem:f-val-d}, and \ref{lem:shadow},
Theorem~\ref{thm:drs-d}
and Lemmas \ref{lem:p-imp-dir}, 
\ref{lem:f-val-d}, and \ref{lem:shadow},
and 
Theorem~\ref{thm:drs-e} and 
Lemma~\ref{lem:f-val-d}, and \ref{lem:shadow3}.
Combining the three results into one
gives us Theorem~\ref{thm:admm2}.
\qed

\section{When strong duality fails}
\label{ss:sd-fail}
In the analyses of DRS,
we assumed strong duality holds.
When strong duality fails, i.e., when $d^\star<p^\star$,
we conjecture that DRS fails.

\begin{conj}
When strong duality fails, DRS necessarily fails in that
\[
\liminf_{k\rightarrow\infty}
f(x^{k+1/2})+g(x^{k+1})< p^\star.
\]
In other words, DRS finds the wrong objective value.
\end{conj}

As discussed in Section~\ref{ss:drs-conv}, DRS tries to achieve feasibility and optimality.
As discussed in Section~\ref{ss:p-subval}, strong duality is well-posedness.
Therefore, when the problem is ill-posed, 
we expect DRS to reduce the function value below $p^\star$
while achieving an infinitesimal infeasibility.
We support the conjecture with examples.

We first present an analytical counter example.
Consider the problem taken from \cite{liu2017new}
\begin{equation*}
\begin{array}{ll}
\mbox{minimize}& 
\underbrace{\delta_{
\{(x_1,x_2,x_3)\,|\,x_3\ge (x_1^2+x_2^2)^{1/2}\}}(x)}_{f(x)}+\underbrace{x_1+\delta_{\{(x_1,x_2,x_3)\,|\,x_2=x_3\}}}_{g(x)}(x)
\end{array}
\end{equation*}
which has the solution set $\{(0,t,t)\,|\,t\in \reals\}$
and optimal value $p^\star=0$.
Its dual problem
\[
\begin{array}{ll}
\mbox{maximize} &
-\delta_{\{(\nu_1,\nu_2,\nu_3)\,|\,-\nu_3\ge (\nu_1^2+\nu_2^2)^{1/2}\}}(\nu)
-\delta_{\{(\nu_1,\nu_2,\nu_3)\,|\,\nu_1=1,\,\nu_2=-\nu_3\}}(-\nu)
\end{array}
\]
is infeasible.
Given $z^0=(z^0_1,z^0_2,0)$, 
the DRS iterates have the form
\begin{align*}
 z^{k+1}_1&=\frac{1}{2}z^k_1-\gamma \\
 z^{k+1}_2&=\frac{1}{2}z^k_2+\frac{1}{2}\sqrt{(z^k_1)^2+(z^k_2)^2}\\
 z^{k+1}_3&=0.
\end{align*}
With this, it is relatively straightforward to show
$x^{k+1/2}-x^{k+1}\rightarrow \vz$,
$x^{k+1/2}_1\rightarrow -2\gamma$, $x^{k+1/2}_2\rightarrow \infty$, $x^{k+1/2}_3\rightarrow \infty$,
and $f(x^{k+1/2})+g(x^{k+1})\rightarrow -2\gamma$.
Also, 
$x^{k+1/2}\nrightarrow \dom f\cap \dom g$ even though $x^{k+1/2}-x^{k+1}\rightarrow \vz$.


Note that 
\[
d^\star<\lim_{k\rightarrow\infty}f(x^{k+1/2})+g(x^{k+1})<p^\star.
\]
So this counterexample proves, at least in some cases, 
that DRS solves neither the primal nor the dual problem
in the absence of strong duality.

Next, we present more experimental counter examples that support the conjecture.
We run DRS on these problems report the experimental results.

The problem, taken from \cite{bertsekas2009},
\[
\underset{x\in \mathbb{R}^2}{\mbox{minimize}}
\quad
\underbrace{
\exp(-\sqrt{x_1x_2})}_{f(x)}
+
\underbrace{\delta_{\{(x_1,x_2)\,|\,x_1=0\}}(x)}_{g(x)}
\]
has  $p^\star=1$ but $d^\star=0$.
Experimentally, 
for all $\gamma>0$ and choice of $z^0$
we observe $d^\star<\lim _{k\rightarrow \infty}f(x^{k+1/2})+g(x^{k+1})<p^\star$.

The problem, taken from \cite{drusvyatskiy2017},
\[
\underset{X\in \mathbf{S}^3}{\mbox{minimize}}
\quad
\underbrace{
\delta_{ \mathbf{S}^3_+}(X)}_{f(X)}+
\underbrace{X_{22}+\delta_{\{
X\in \mathbf{S}^3 \,|\, X_{33} = 0, X_{22} + 2X_{13} = 1
\}}(X)}_{g(X)},
\]
where $\mathbf{S}^3$ and $\mathbf{S}^3_+$ respectively denote
the set of symmetric and positive semidefinite
$3\times 3$ matrices,
has  $p^\star=1$ but $d^\star=0$.
Experimentally, we observe
$d^\star=\lim_{k\rightarrow \infty }f(x^{k+1/2})+g(x^{k+1})$ for $\gamma\geq 0.5$,
and  $d^\star<\lim_{k\rightarrow \infty }f(x^{k+1/2})+g(x^{k+1})<p^\star$ for $0<\gamma< 0.5$.
This behaviour does not depend on $z^0$.

The problem, taken from \cite{ye2004linear},
\[
\underset{X\in \mathbf{S}^3}{\mbox{minimize}}
\quad
\underbrace{
\delta_{ \mathbf{S}^3_+}(X)}_{f(X)}+
\underbrace{2X_{12}+\delta_{\{
X\in \mathbf{S}^3 \,|\, X_{22} = 0, -2X_{12}+2X_{33}=2
\}}(X)}_{g(X)}
\]
has  $p^\star=0$ but $d^\star=-2$.
Experimentally, we observe $d^\star=\newline\lim_{k\rightarrow \infty }f(x^{k+1/2})+g(x^{k+1})$ for $\gamma\geq 1$, and  $d^\star<\lim_{k\rightarrow \infty }f(x^{k+1/2})+g(x^{k+1})<p^\star$ for $0<\gamma< 1$. This behaviour does not depend on $z^0$.

The problem, taken from \cite{tunccel2012strong},
\[
\underset{X\in \mathbf{S}^5}{\mbox{minimize}}
\quad
\underbrace{
\delta_{ \mathbf{S}^5_+}(X)}_{f(X)}+
\underbrace{X_{44}+X_{55}+\delta_{\{
X\in \mathbf{S}^3 \,|\, X_{11} = 0, X_{22}=1, X_{34}=1, 2X_{13}+2X_{45}+X_{55}=1
\}}(X)}_{g(X)}
\]
has  $p^\star=(\sqrt{5}-1)/2$ but $d^\star=0$.
Experimentally, we observe $d^\star=\newline\lim_{k\rightarrow \infty }f(x^{k+1/2})+g(x^{k+1})$ for $\gamma\geq 0.8$, and  $d^\star<\lim_{k\rightarrow \infty }f(x^{k+1/2})+g(x^{k+1})<p^\star$ for $0<\gamma< 0.8$. This behaviour does not depend on $z^0$.

The conjecture holds for all examples.
Interestingly, for some examples, there is a 
threshold $\gamma_\mathrm{min}$
such that
$d^\star<\lim_{k\rightarrow \infty }f(x^{k+1/2})+g(x^{k+1})<p^\star$
when $0<\gamma < \gamma_\mathrm{min}$
and
$d^\star=\lim_{k\rightarrow \infty }f(x^{k+1/2})+g(x^{k+1})$
when $\gamma_\mathrm{min}\le \gamma$.
We do not have an explanation for this phenomenon.

\section{Conclusion}
\label{s:conclusion}
In this paper, we analyzed DRS and ADMM under pathologies.
We show that when strong duality holds, the iterates of DRS and ADMM are approximately feasible and approximately optimal in the sense discussed in Sections \ref{ss:interpretation} and \ref{ss:admm-interpretation}.
Furthermore, we conjectured that DRS necessarily fails when strong duality fails, and we provided empirical evidence supporting this conjecture.

As discussed in Section~\ref{ss:sd-fail},
DRS exhibits an interesting behavior in the absence of strong duality,
and we do not have an explanation for it.
Analyzing this phenomenon and addressing the conjecture
is an interesting direction of future research.

%

%

For non-pathological problems, DRS can be generalized  
with an over-under relaxation parameter between $0$ and $2$.
The pathological DRS analysis of this paper immediately extends to this generalized setup.
For non-pathological problems, ADMM can be generalized 
with an over-under relaxation parameter between $0$ and $(1+\sqrt{5})/2$.
This generalization arises when ADMM is analyzed directly through a Lyapunov function, 
and not through DRS
\cite{fortin1983,glowinski1984,bertsekaspar,fezel2013,Deng2016,Chen2017,chen2018}.
The pathological ADMM analysis of this paper \emph{does not} immediately extend to this generalized setup.
Analyzing this form of ADMM applied to pathological problems is an interesting direction of future research.

\begin{acknowledgements}
This work was partially supported by NSF DMS-1720237 and ONR N0001417121.
\end{acknowledgements}

\bibliographystyle{spmpsci}      
\bibliography{./drs}   

\end{document}